\documentclass[12pt]{article}
\usepackage{amsmath}
\usepackage{amsthm}
\usepackage{amssymb}
\usepackage{latexsym}
\usepackage{amsfonts}
\newcommand{\nc}{\newcommand*}
\newcommand{\rnc}{\renewcommand*}
\nc{\ket}[1]{{\vert{#1}\rangle}}                
\nc{\bra}[1]{{\langle{#1}\vert}}                
\nc{\braket}[2]{{\langle{#1}\vert{#2}\rangle}}  
\nc{\ketbra}[2]{{\vert{#1}\rangle\langle{#2}\vert}}  
\nc{\kett}[2]{{\vert{#1,#2}\rangle}}            
\nc{\braa}[2]{{\langle{#1,#2}\vert}}            
\nc{\spr}[2]{{\left({#1},{#2}\right)}}          
\nc{\lcr}[2]{{\left[{#1},{#2}\right]}}          
\nc{\ketx}[2]{{\vert{#1}\rangle}_{#2}}          
\nc{\brax}[2]{{\langle{#1}\vert}_{#2}}          
\nc{\braketx}[3]{{{}_{#3}\langle{#1}\vert{#2}\rangle_{#3}}}  
\nc{\brkt}[3]{{\langle{#1}\vert{#3}\vert{#2}\rangle}}  
\nc{\ketbrax}[3]{{\vert{#1}\rangle_{#3}{}_{#3}\langle{#2}\vert}} 
\nc{\ketbray}[3]{{\vert{#1}\rangle{\rule[-5pt]{0pt}{1pt}}_{#3}{\rule[-5pt]{0pt}{1pt}}_{#3}\langle{#2}\vert}} 
\nc{\skob}[1]{{\left({#1}\right)}}              
\nc{\bskob}[1]{{\left[{#1}\right]}}             
\nc{\cskob}[1]{{\left\{{#1}\right\}}}           
\nc{\modx}[1]{{\vert{#1}\vert}}                 
\nc{\MO}[1]{{\langle{#1}\rangle}}               
\nc{\CKO}[1]{{\bskob{\langle{#1}^2\rangle - {\langle{#1}\rangle}^2}}^{\half}}       
\nc{\bmdx}[1]{{\left|{#1}\right}}                 
\nc{\derx}[1]{\frac{\rm d}{{\rm d}\,{#1}}}         
\nc{\deryx}[2]{\frac{{\rm d}{#1}}{{\rm d}\,{#2}}}   
\nc{\dertx}[1]{\frac{{\rm d}^2}{{\rm d}\,{#1}^2}}         
\nc{\dertyx}[2]{\frac{{\rm d}^2\,{#1}}{{\rm d}\,{#2}^2}}   
\nc{\dern}[2]{\frac{{\rm d}^{#1}}{{\rm d}\,{#2}^{#1}}}     
\nc{\ddz}[1]{{\rm d}(\text{Re}{#1}){\rm d}(\text{Im}{#1})}
\nc{\cald}{{\cal D}} \nc{\calf}{{\cal F}} \nc{\calh}{{\cal H}}
\nc{\calk}{{\cal K}} \nc{\caln}{{\cal N}} \nc{\Gz}{{G^{(0)}}}
\nc{\Go}{{G^{(1)}}} \nc{\ebf}{{\bf E}} \nc{\dbf}{{\bf D}}
\nc{\kbf}{{\bf K}} \nc{\nbf}{{\bf N}} \nc{\cals}{{\cal S}}
\nc{\calx}{{\cal X}} \nc{\calw}{{\cal W}} \nc{\mI}{{\mathbb{I}}}
\nc{\mN}{{\mathbb{N}}} \nc{\mR}{{\mathbb{R}}}
\nc{\mC}{{\mathbb{C}}} \nc{\mS}{{\mathbb{S}}} \nc{\bone}{\,{\rm
I}\hspace{-7pt}1} \nc{\gmf}[1]{\varGamma({ #1})}
\nc{\Szi}[1]{\sum_{ #1 = 0 }^{\infty}} \nc{\Sn}{\Szi{n}}
\nc{\Sm}{\Szi{m}} \nc{\Szx}[2]{\sum_{ #1 = 0 }^{#2}}
\nc{\Syx}[3]{\sum_{ #1 = #2 }^{#3}} \nc{\Sud}[2]{\sum^{#1}_{#2}}
\nc{\Fto}[4]{\raisebox{-3pt}{${}_2\text{\large
F}_1$}\left(\ts{\genfrac{}{}{0pt}{}{#1, #2}{#3}\left.\right| #4
}\right)} \nc{\BFto}[4]{\raisebox{-3pt}{${}_2\text{\large
F}_1$}\left(\ts{\genfrac{}{}{0pt}{}{#1, #2}{#3}\biggl.\biggr| #4
}\right)} \nc{\aFto}{\Fto{\half}{\frac32}{1}{2J}}
\nc{\bFto}{\Fto{\frac32}{\frac52}{1}{2J}}
\nc{\dFto}{\Fto{\frac32}{\frac52}{2}{2J}}
\nc{\xFto}[1]{\raisebox{-3pt}{${}_2\text{\Large F}_1$}\left(
\ts{\genfrac{}{}{0pt}{}{\half, \frac32}{1}\biggl|\biggr. #1
}\right)} \nc{\Fzo}[2]{\ts{\raisebox{-3pt}{${}_0\text{\Large
F}_1$}\left( \genfrac{}{}{0pt}{}{-}{#1}\Biggl|\Biggr. #2 \right)}}
\nc{\Ff}[6]{\raisebox{-3pt}{$\text{\large F}_4$}\left(
#1,#2;#3,#4;#5,#6\right)}
\nc{\Ftf}[7]{\raisebox{-3pt}{${}_2\text{\large
F}_4$}\left(\ts{\genfrac{}{}{0pt}{}{#1,
#2}{#3,#4,#5,#6}\left|\right. #7 }\right)} \nc{\Izi}[2]{\int_{ 0
}^{\infty}{#1}{\rm d}\,{#2}} \nc{\Izx}[3]{\int_{ 0 }^{#1}{#2}{\rm
d}\,{#3}} \nc{\Ixy}[4]{\int_{#1}^{#2}{#3}{\rm d}\,{#4}}
\nc{\Ioo}[2]{\int_{-1}^{1}{#1}{\rm d}\,{#2}}
\nc{\ophi}[1]{\varphi^{(0)}_{#1}}
\nc{\opsi}[1]{\psi^{\text{opt}}_{#1}}
\nc{\gexp}[2]{\exp_{#1}\left({#2}\right)}          
\nc{\eby}[1]{\stackrel{\reff{#1}}{=}} \nc\dagx[1]{{#1}^\dagger}
\rnc{\Re}{\text{Re}} \rnc{\Im}{\text{Im}}
\def\reff#1{(\ref{#1})}
\def\ds{\displaystyle}
\def\ts{\textstyle}
\def\half{{\frac{1}{2}}}

\def\ra{\rightarrow}
\def\te{\text{e}}

\nc{\ocgk}{\int\ketbra{J,\gamma}{J,\gamma}{\rm d}\,\mu(J,\gamma)}
\nc{\gm}[2]{\gamma_{#1,#2}}
\oddsidemargin=0.25cm \rnc{\baselinestretch}{1.2}
\begin{document}
$${}$$
\begin{flushright}
{Translation from Russian}\\ {of the article published in}\\ {{\it
ZNS POMI}, {\bf 285}, 39-52 (2002)}
\end{flushright}
$${}$$

\centerline{\huge\bf Coherent states}

\medskip

\centerline{\large\bf for the}
\bigskip

\centerline{{\huge\bf Legendre oscillator} \large\footnote{This
research is supported in part by RFFI grant no. 00-01-00500.}}

\bigskip
\smallskip
{\flushleft{\large\bf V.V.Borzov,${}^*$
E.V.Damaskinsky${}^{**}$}}
\smallskip
{\flushleft{${}^*$ St.-Petersburg University of
Telecommunications\\ E-mail: vadim@VB6384.spb.edu}} \vspace{-.1cm}
{\flushleft{${}^{**}$St.-Petersburg University of the Defence
Engineering Constructions\\ E-mail: evd@pdmi.ras.ru}}
\bigskip

\centerline{\bf Abstract}
\begin{quote}
A new oscillator-like system called by the Legendre oscillator is
introduced in this note.
 The two families of coherent states (coherent states as
eigenvectors of the annihilation operator and the Klauder ---
Gazeau temporally stable coherent states) are defined and
investigated for this oscillator.
\end{quote}
\vskip 0.5cm

\section{Introduction}
In the present note we  construct  coherent states for an
oscillator-like system called by the Legendre oscillator.Note that
the Legendre polynomials play the same role for this system as the
Hermite polynomials for standard boson oscillator.

It is known that in the case of harmonic oscillator three standard
definitions of coherent state (as eigenvector of an annihilation
operator; as state generated by a shift operator from the vacuum
state and as state minimizing the uncertainty relation) are
equivalent.This means that they generate the same set of states.
However, in the general case this is not so. Below we shall
construct coherent state connected with the Legendre polynomials.
This construction
results from the analysis~\cite{borz} of the linkage between
orthogonal polynomials and generalized oscillator
algebras~\cite{damkul,bordamkul}. Namely, we shall define for the
Legendre oscillator analogues of Barut - Girardello coherent
states~\cite{BarutGirardello} and Gazeau - Klauder coherent
states~\cite {GazeauKlauder99a}.
%
The more detailed exposition is demanded for definition of
coherent states of the Perelomov-type
~\cite{PerelomovCMP,PerelemovUFN, PerelomovBk}. So it will be
postponed to the other publication.

\section{The Legendre oscillator.}
For the reader convenience we  remind some information about the
Legendre polynomials.
 The Legendre polynomials $P_{(x)}$ are the
solutions of the differential equations
\begin {equation}
(1-x^2)y^{\prime\prime}-2xy^{\prime}+n(n+1)y=0,\qquad
y(x)=P_{n}(x)
\end {equation}
which satisfy the orthogonality condition
\begin{equation}
\Ioo{P_n(x)P_m(x)}{x}=\frac{2}{2n+1}\delta_{n\;m}. \label{le2a}
\end{equation}
They also are solutions of the following recurrent relations
\begin{equation}
(2n+1)xP_n(x)=(n+1)P_{n+1}(x)+nP_{n-1}(x);\qquad P_{0}(x)=1;
\qquad (P_{-1}(x)=0), \qquad n\in\mN_0. \label{le2}
\end{equation}
The Legendre polynomials are defined by the relation ($ n\in\mN _
0 $)
\begin{align}
P_n(x)=\Fto{-n}{\,n+1}{1}{\frac{1-x}2}&=\Szx{m}{[\![\frac{n}{2}]\!]}
       \frac{(-1)^m}{2^n}\binom{n}{m}\binom{2n-2m}{n}x^{n-2m}=\nonumber\\
&=\Szx{m}{[\![\frac{n}{2}]\!]}
       \frac{(-1)^m(2n-2m)!}{2^nm!(n-m)!(n-2m)!}x^{n-2m},
\label{le1}
\end{align}
where a symbol $[\![n]\!]$ denotes the integer part of a number
$n.$ Below we shall use the following generating function for this
polynomials
\begin{equation}\label{le2b}
\Szi{n}\frac{(\gamma)_n}{n!}P_n(x)z^n=(1-xz)^{-\gamma}
\Fto{\half\gamma}{\,\half(1+\gamma)}{1}{\frac{(x^2-1)z^2}{(1-xz)^2}};
\quad |x|\leq 1,\, |z|\leq 1.
\end{equation}

Simultaneously with the Legendre polynomials we  consider the
Legendre functions
\begin{equation}
\psi_n(x)=\sqrt{2n+1}P_{n}(x),\qquad n\in\mN_0 , \label{le3a}
\end{equation}
which form a orthonormal basis $\left\{
\ket{n}\equiv\psi_{n}(x)\right\}_{n=0}^{\infty}$ in the Hilbert
space
\begin{equation}
\calh :=L^2\left( [-1,1] , {\ts\half}\mathrm{d\,}x\right).
\label{le4}
\end{equation}
These functions fulfill the recurrent relations
\begin{equation}
x\psi_n(x)=b_{n-1}\psi_{n-1}(x)+b_{n}\psi_{n+1}(x),\qquad
\psi_{-1}(x)=0,\quad\psi_{0}(x)=1, \label{le5}
\end{equation}
with coefficients
\begin{equation}
b_n=\sqrt{\ds\frac{(n+1)^2}{(2n+1)(2n+3)}},\qquad n\geq 0 .
\label{le6}
\end{equation}

In the given research the Legendre polynomials  $P_{n}(x)$ and the
Legendre functions $\psi_n(x)$ play the same role as the Hermite
polynomials and the Hermite functions play in the standard quantum
mechanics.

In the Hilbert space $\calh$ we  define the generalized position
operator $X$ connected with the Legendre polynomials
 $P_{n}(x)$ as an operator of multiplication by argument:
\begin{equation}
X\ket{n}=x\ket{n}. \label{le8}
\end{equation}
 Taking into account a relation \reff{le5},we have
\begin{equation}
X\psi_{n}(x)= b_{n}\psi_{n+1}(x)+b_{n-1}\psi_{n-1}(x), \label{le9}
\end{equation}
where the coefficients $b_n$ are defined by the relation
\reff{le6}. Because $\sum_{k=0}^\infty \frac 1{b_k}=\infty ,$ the
operator $X$ is a selfajoint operator in the space $\calh$ (see
\cite {Akhieser, BirmanS, bordamkul}) .

Let us define a generalized momentum operator $P$ by the way
described in~\cite{borz}. The operator $P$ acts on the basis
elements in $\calh$  by the following formula
\begin{equation}
P\ket{n}=i\skob{b_{n}\ket{n+1}-b_{n-1}\ket{n-1}}. \label{le12}
\end{equation}

Calculating usual commutator of operators $X$ and $P$ on the basis
elements
 , we obtain
\begin{equation}\label{le14}
\lcr{X}{P}\ket{n}=2i\skob{{b_{n}}^2-{b_{n-1}}^2}\ket{n}=
\frac{2i}{(2n-1)(2n+1)(2n+3)}\ket{n} .
\end{equation}

Now we define the creation and annihilation operators by the
standard relations
\begin{equation}\label{le15}
a^{(+)}=\frac{1}{\sqrt{2}}\skob{X-iP},\qquad
a^{(-)}=\frac{1}{\sqrt{2}}\skob{X+iP} .
\end{equation}
On the basis elements in $ \calh $ these operators act by the rule
\begin{equation}\label{le16}
a^{(+)}\ket{n}=\sqrt{2}b_{n}\ket{n+1},\qquad
a^{(-)}\ket{n}=\sqrt{2}b_{n-1}\ket{n-1}.
\end{equation}
They satisfy the commutation relations
\begin{equation}\label{le17}
\lcr{a^{(-)}}{a^{(+)}}= \frac1i\lcr{X}{P} .
\end{equation}

Now we introduce the state numbering operator $N$ and Hamiltonian
$H$, by the following formulae
\begin{equation}\label{le18}
N\ket{n}=n\ket{n},\qquad H=X^2+P^2=a^{(+)}a^{(-)}+a^{(-)}a^{(+)} .
\end{equation}
The eigenvalues of the operator $H$ are equal to
\begin{equation}\label{le18 }
\lambda_0=2{b_0}^2,\qquad \lambda_n=2\skob{{b_{n-1}}^2+{b_{n}}^2}.
\end{equation}

It is natural call the introduced system by the Legendre
oscillator.

\section{Barut - Girardello coherent states for the Legendre oscillator}
In this section we define the coherent states for the Legendre
oscillator in the space $\calh$ as eigenvectors of the
annihilation operator $a^{(-)}$
\begin{equation}
a^{(-)}\ket{z}=z\ket{z}. \label{le19}
\end{equation}
It is known that
\begin{equation}
\ket{z}=\mathcal{N}^{-1}\Szi{n}
\frac{z^n}{\left(\sqrt{2}b_{n-1}\right)!}\ket{n}. \label{le14a}
\end{equation}
The normalizing factor is equal to
\begin{equation}
\qquad \mathcal{N}^2=\braket{z}{z}=
\Szi{n}\frac{|z|^{2n}}{\left(2{b_{n-1}}^2\right)!}\equiv
\gexp{[2{b_{n-1}}^2]}{|z|^{2}} . \label{le20}
\end{equation}
Because of\footnote{Here $(a)_n$ is the Pochhammer symbol defined
by the relation $$ (a)_0=1, (a)_n=a(a+1)\cdots
(a+n-1)=\frac{\gmf{a+n}}{\gmf{a}}, \quad n=1,2,\ldots .$$}
\begin{equation}
\left(2{b_{n-1}}^2\right)!=\frac{2^n(n!)^2}{(2n-1)!!(2n+1)!!}
=\frac{(n!)(1)_n}{2^n(\half)_n(\frac32)_n}, \label{Le.16}
\end{equation}
the radius of convergence of a series \reff{le20} equals to
$\ds\frac{1}{\sqrt{2}}$ and
\begin{equation}
\mathcal{N}^2=
\Szi{n}\frac{(\half)_n(\frac32)_n}{n!(1)_n}(2|z|^{2})^n=
\BFto{\half}{\frac32}{1}{2|z|^{2}}. \label{Le.17}
\end{equation}
Substituting \reff{Le.17} in \reff{le14a} and using a relation
$\ket{n}=\psi_n(x)=\sqrt{2n+1}P_n(x),$ we obtain
\begin{align}
\ket{z}&=\left[ \BFto{\half}{\frac32}{1}{2|z|^{2}}
\right]^{-\half}
\Szi{n}\frac{(2n-1)!!(2n+1)}{n!}P_n(x)z^n=\nonumber \\ &=\left[
\BFto{\half}{\frac32}{1}{2|z|^{2}} \right]^{-\half}
\Szi{n}\frac{(\frac{3}{2})_n(\sqrt{2}z)^n}{n!}P_n(x).
\label{Le.18}
\end{align}
From\reff{le2b} as $\gamma=\frac32$ and $z\rightarrow 2z$ we
obtain
\begin{equation}
\ket{z}=\left[
\BFto{\half}{\frac32}{1}{2|z|^{2}}\right]^{-\half}\,
\BFto{\frac{3}{4}}{\frac{5}{4}}{1}{\ds\frac{(x^2-1)2z^{2}}{(1-\sqrt{2}xz)^2}}
\,(1-\sqrt{2}xz)^{-\frac{3}{2}} . \label{Le.19}
\end{equation}

Our following task is to construct a measure $${\rm
d}\mu(|z|^2)=W(|z|^2){\rm d}^2z, \quad\text{such that}\quad
\int\!\!\int_{\mathbb{C}}W(|z|^2)\,\ketbra{z}{z}\,{\rm d}^2z=1, $$
where
 ${\rm d}^2z={\rm d}(\Re{z}){\rm d}(\Im{z}).$
It is known (see, for example, \cite{SixdenierrsPensonSolomon}),
that this problem is reduced to a solution of the following the
Hausdorf moment problem
\begin{equation}
\Szi{n}\frac{\pi}{(2{b_{n-1}}^2)!}\left[\Izx{\half}{t^nW(t)}{t}\right]
\ketbra{n}{n}=1,\qquad (t=|z|^2) \label{Le.20}
\end{equation}
or
\begin{equation}
\Izx{\half}{t^nW(t)}{t}=\frac{1}{\pi}(2{b_{n-1}}^2)!.
\label{Le.21}
\end{equation}
Substituting \reff{Le.16} in \reff{Le.21}, we have (for $\tau=2t$)
$$ \half\Izx{1}{{\tau}^nW(\ts{\half}\tau)}{\tau}=
\frac{1}{\pi}\frac{(n!)^2}{(\half)_n(\frac32)_n} . $$
So, it is
necessary to solve a following Hausdorf moment problem
\begin{equation}
\half\Izx{1}{{\tau}^nW(\ts{\half}\tau)}{\tau}=
\frac{1}{\pi}\frac{\left(\gmf{n+1}\right)^2\gmf{\half}\gmf{\frac32}}
{\gmf{n+\half}\gmf{n+\frac32}} \label{Le.22}
\end{equation}
or, taking into the account $\gmf{\half}\gmf{\frac32}=\half\pi$,
\begin{equation}
\Izx{1}{{\tau}^nW(\ts{\half}\tau)}{\tau}=
\frac{\left(\gmf{n+1}\right)^2} {\gmf{n+\half}\gmf{n+\frac32}},
\qquad n\geq 0 . \label{Le.23}
\end{equation}

Making in an integral (see (7.127) in \cite{G-R})
\begin{equation}
J=\Ixy{-1}{1}{(1+x)^{\sigma}P_{\nu}(x)}{x}=
\frac{2^{1+\sigma}\left(\gmf{1+\sigma}\right)^2}
{\gmf{2+\sigma+\nu}\gmf{1+\sigma-\nu}}, \qquad \Re{\sigma}>-1,
\label{Le.24}
\end{equation}
the replacement $x=2\tau-1$, we  receive
$$J=\Izx{1}{2^{\sigma+1}{\tau}^{\sigma}P_{\nu}(2\tau-1)}{\tau}.$$
This allows us to rewrite \reff{Le.24} in the form
\begin{equation}
\Izx{1}{{\tau}^{\sigma}P_{\nu}(2\tau-1)}{\tau}=
\frac{\left(\gmf{\sigma+1}\right)^2}
{\gmf{\sigma+\nu+2}\gmf{1+\sigma-\nu}}. \label{Le.24a}
\end{equation}
Choosing $\sigma=n$ AND $\nu=\half,$  we obtains
\begin{equation}
\Izx{1}{{\tau}^{n}\left(P_{\half}(2\tau-1)-{\tau}^{\half}\right)}{\tau}
+\Izx{1}{{\tau}^{n+\half}}{\tau}= \frac{\left(\gmf{n+1}\right)^2}
{(n+\frac32)\gmf{n+\half}\gmf{n+\frac32}}.
\label{Le.25}\end{equation} Let's denote
\begin{equation}
\tau^{-\half}P_{\half}(2\tau-1)-1= \Ixy{\tau}{1}{q(t)}{t}.
\label{Le.26}
\end{equation}
After differentiating, we find
function $q(\tau)$ for $0<\tau<1$
\begin{equation}
q(\tau)=-\left({\tau}^{-\half}P_{\half}(2\tau-1)\right)^{\prime}.
\label{Le.27}
\end{equation}

For an investigation of a singularity, arising at $\tau\rightarrow
0^+$, we shall consider an integral
$$\Izx{1}{{\tau}^{n+\half}\left(\Ixy{\tau}{1}{q(t)}{t}
\right)}{\tau}.$$
Integrating by parts, we obtain
\begin{equation}
\Izx{1}{{\tau}^{n+\half}\left(\Ixy{\tau}{1}{q(t)}{t}
\right)}{\tau}=
\frac{{\tau}^{n+\frac32}}{n+\frac32}\Ixy{\tau}{1}{q(t)}{t}
\Bigl|^1_{0}\Bigr.
+\Izx{1}{q(\tau)\frac{{\tau}^{n+\frac32}}{n+\frac32}}{\tau}.
\label{Le.28}\end{equation} Let's remark, that $$
\Ixy{\tau}{1}{q(t)}{t} \xrightarrow[\tau\rightarrow 1]{}0
\qquad\text{and}\qquad \tau^{3/2}\Ixy{\tau}{1}{q(t)}{t}
\xrightarrow[\tau\rightarrow 0]{}0, $$ from which it follows, that
the term outside the integral in \reff{Le.28} is equal to zero.
From \reff{Le.25} and \reff{Le.26} it follows, that the integral
in the left hand side of the relation \reff{Le.28} is equal to $$
\Izx{1}{{\tau}^{n+\half}\left(\Ixy{\tau}{1}{q(t)}{t}
\right)}{\tau} +\frac{1}{n+\frac32}=
\frac{\left(\gmf{n+1}\right)^2}
{(n+\frac32)\gmf{n+\half}\gmf{n+\frac32}}, $$ so that for $n\geq
0$ we have
\begin{equation}
\Izx{1}{q(\tau){\tau}^{n+\frac32}}{\tau}+1=
\frac{\left(\gmf{n+1}\right)^2} {\gmf{n+\half}\gmf{n+\frac32}}.
\label{Le.29}
\end{equation}
Thus the moment problem is solved by the distribution
$$
W(\ts{\half}\tau)={\tau}^{\frac32}q(\tau)+2\delta(\tau-1)=
-{\tau}^{\frac32}\left({\tau}^{-\half}P_{\half}(2\tau-1)\right)^{\prime}
+2\delta(\tau-1) $$ or
\begin{equation}
W(t)=-(2t)^{\frac32} \left(
(2t)^{-\half}P_{\half}(4t-1)\right)^{\prime} +2\delta(2t-1),
\qquad 0<t\leq\half . \label{Le.30a}
\end{equation}
Using the formula (8.832 (1)) from \cite{G-R}, we obtain
\begin{equation}
W(t)=\frac{(16t-5)P_{\half}(4t-1)-3P_{\frac32}(4t-1)}{2(2t-1)}+2\delta(2t-1),
\qquad 0<t\leq\half . \label{Le.30}\end{equation} Finally, for
$0<|z|\leq\frac{1}{\sqrt{2}}$ we have
\begin{equation}
{\rm d}\mu(|z|^2)=\left[
\frac{(16|z|^2-5)P_{\half}(4|z|^2-1)-3P_{\frac32}(4|z|^2-1)}{2(2|z|^2-1)}
+2\delta(2|z|^2-1)\right] \ddz{z} . \label{Le.31}\end{equation}

We calculate overlap of two coherent states
\begin{align}\label{}
\braket{z_1}{z_2}&=
\bskob{\Fto{\half}{\frac32}{1}{2\modx{z_1}^2}\,
\Fto{\half}{\frac32}{1}{2\modx{z_2}^2}}^{-\half}
\Szi{n}\frac{{\bar{z}_1}^n{z_2}^n}{\skob{2{b_{n-1}}^2}!}\eby{Le.16}
\nonumber \\ &=\bskob{\Fto{\half}{\frac32}{1}{2\modx{z_1}^2}\,
\Fto{\half}{\frac32}{1}{2\modx{z_2}^2}}^{-\half}
\Szi{n}\frac{\skob{2\bar{z}_1z_2}^n\skob{\half}_n\skob{\frac32}_n}
{(n!)(1)_n}\eby{Le.17} \nonumber \\
&=\bskob{\Fto{\half}{\frac32}{1}{2\modx{z_1}^2}\,
\Fto{\half}{\frac32}{1}{2\modx{z_2}^2}}^{-\half}
\Fto{\half}{\frac32}{1}{2\bar{z}_1z_2} . \label{Le.32}\end{align}

To arbitrary normalized state
$\ket{f}=\Szi{n}f_n\ket{\psi_n}\in\calh $ ($\Szi{n}{f_n}^2=1$) we
can put in correspondence a function analytical on
${\mC}_{1/\sqrt{2}}$, by the rule
\begin{equation}\label{Le.34}
f(z)=\caln(z)\braket{z}{f}=
\Szi{n}\sqrt{\ts{\left(\half\right)_n\,
\left(\frac32\right)_n}}\,\,\frac{f_n}{n!}\,(2z)^n ,
\end{equation}
so that after expansion on coherent states we have
\begin{equation}\label{Le.35}
\ket{f}=\int_{{\mC}_{1/\sqrt{2}}}\braket{z}{f}\,\ket{z} {\rm
d}\mu(|z|^2)=\int_{{\mC}_{1/\sqrt{2}}}
\bskob{\Fto{\half}{\frac32}{1}{2\modx{z}^2}}^{-\half}f(z)\,
\ket{z}{\rm d}\mu(|z|^2)
\end{equation}
and
\begin{equation}\label{Le.36}
\braket{f}{f}= \int_{{\mC}_{1/\sqrt{2}}}
\bskob{\Fto{\half}{\frac32}{1}{2\modx{z}^2}}^{-1}\modx{f(z)}^2
{\rm d}\mu(|z|^2)<\infty .
\end{equation}

\section{Klauder - Gazeau coherent states for the Legendre oscillator.}
The Klauder - Gazeau temporary stable coherent states~\cite
{GazeauKlauder99a} is convenient to apply in a case when the
Hamiltonian $H$ is nonlinear and its terms are not generators of
group of a symmetry. These states can be defined by the relation
\begin{equation}\label{Le.38}
\ket{J, \gamma}:=\caln(J)^{-1}\Szi{n}\frac{J^{\frac
n2}}{\sqrt{\rho_n}} \te^{-i\gamma\lambda_n}\ket{\psi_n},
\end{equation}
where $\rho_n=\lambda_{1}\lambda_{2}\cdot\ldots\cdot\lambda_{n},
\quad n\geq 1; \quad \rho_0=1,$ and the normalizing coefficient is
equal
\begin{equation}\label{Le.40}
\caln(J)^{2}=\Szi{n}\frac{J^{n}}{\rho_n}.
\end{equation}
Parameters $J$ and $\gamma$ takes the values $J\geq 0,$
$\gamma\in\mR$.Note that these parameters are generalization of
the module and argument (extended up to an infinite covering of a
segment $\bskob{0;2\pi}$)of a standard parameter
$z=\modx{z}\te^{i\gamma}$ of coherent states. One can to consider
these parameters as analogue of the classical action - angle
variables.

In  considered case the Hamiltonian $H$  is a positive selfajoint
operator in a Hilbert space $\calh$ with a simple discrete
(ordered by decreasing) spectrum
 $\left\{ \lambda_n=2{b_{n-1}}^2\right\}_{n=1}^{\infty}$,
where $b_n$ is taken from \reff{le6}, so that we can write
\begin{gather}
\rho_0=1,\, \rho_n=\skob{2{b_{n-1}}^2}!=
\frac{(n!)(1)_n}{2^n\skob{\half}_n\skob{\frac32}_n},\quad n\geq 1
, \label{Le.41}\\
\caln(J)^{2}=\Szi{n}J^{n}\frac{2^n\skob{\half}_n\skob{\frac32}_n}{(n!)(1)_n}=
\Fto{\half}{\frac32}{1}{2J} . \label{Le.42}
\end{gather}
The radius of convergence of a series in \reff{Le.42} is equal to
$R=\lim_{n\ra \infty}\sqrt[n]{\rho_n}=\half.$ Taking into account
\reff{Le.41}, \reff{Le.42}, we obtain
\begin{equation}\label{Le.44}
\ket{J,\gamma}=\frac{1}{\sqrt{\aFto}}\Szi{n}
\frac{\sqrt{2n+1}}{n!} J^{n/2}2^{n/2}
\sqrt{\skob{\half}_n\skob{\frac32}_n}
\te^{-i\sqrt{2}b_{n-1}\gamma}P_n(x),
\end{equation}
or, because  $\skob{\frac32}_n=(2n+1)\skob{\half}_n$,
\begin{equation}\label{Le.45}
\ket{J,\gamma}=\frac{1}{\sqrt{\aFto}}\Szi{n}
\frac{2n+1}{n!}(2J)^{n/2}\skob{\half}_n
\te^{-i\sqrt{2}b_{n-1}\gamma}P_n(x) .
\end{equation}
To proof the validity of the resolution of identity we consider
the relation
\begin{equation}\label{Le.46}
\ocgk=\lim_{T\ra\infty}\frac{1}{2T}\int_T^T{\rm d}\,\gamma
\bskob{\int_0^{\infty}k(J)\ketbra{J,\gamma}{J,\gamma}{\rm d}\,J},
\end{equation}
where
\begin{equation}\label{Le.47}
k(J):=\left\{
\begin{matrix}
{\caln(J)}^2\rho(J) & 0\leq J\leq \half \\
                  0 & J>\half
\end{matrix}
\right. .
\end{equation}
Calculating  an integral over $\gamma$ we obtain
\begin{equation}\label{Le.48}
\ocgk=\Szi{n}\frac{1}{\rho_n}\,\Ixy{0}{\half}{\ketbra{n}{n}}{J}.
\end{equation}
Thus, the resolution of identity
\begin {equation} \label {Le.49}
\ocgk = \bone
\end {equation}
is fulfilled, if the weight function $\rho(J)$ gives a solution of
a moment problem
\begin{equation}\label{Le.50}
\Ixy{0}{\half}{J^n\rho(J)}{J}=\rho_n=
\frac{n!(1)_n}{2^n\skob{\half}_n\skob{\frac32}_n}, \qquad n\geq 0.
\end{equation}
In view of the relations \reff{Le.21} and \reff{Le.30} a solution
of this problem is given by the relation
 ($0<J<\half$, $\gamma\in\mR$)
\begin{equation}\label{Le.51}
\rho(J)=\frac{\pi}{4(2J-1)}\bskob{(16J-5)P_{\half}(4J-1)-3P_{\frac32}(4J-1)}
+\pi\delta(2J-1).
\end{equation}
The temporal stability is obvious, as
\begin{equation}\label{Le.52}
\te^{-iHt}\ket{J,\gamma}=
\frac{1}{\caln}\Szi{n}\frac{J^{n/2}}{\sqrt{\rho_n}}
\te^{-i\gamma\lambda_n}\te^{-it\lambda_n}\ket{\psi_n}=
\ket{J,\gamma+t}.
\end{equation}

We have also
\begin{equation}\label{Le.53}
\brkt{J,\gamma}{J,\gamma}{H} = {\caln}^{-2}(J)\Szi{n}
\frac{{2b_{n-1}}^2}{({2b_{n-1}}^2)!}J^n=J .
\end{equation}

The overlap of two states is given by a relation
\begin{equation}\label{Le.54}
\braket{J^{\prime},\gamma^{\,\prime}}{J,\gamma}=
\frac{1}{N(J)N(J^{\prime})}\Szi{n}\frac{(JJ^{\prime})^{n/2}}{\rho_n}
\exp{\left[-i\lambda_n(\gamma -\gamma^{\prime})\right]} ,
\end{equation}
which, in the concrete case $\gamma =\gamma^{\prime},$ is easily
summarized
\begin{equation}\label{Le.55}
\braket{J^{\prime},\gamma}{J,\gamma}=
\frac{\xFto{2\sqrt{JJ^{\prime}}}} {\sqrt{\aFto\xFto{2J^{\prime}}}}.
\end{equation}

Taking into account a possible physical applications, we
calculate, for  example, some quantities having an immediate
physical sense. So, for an average number  of excitation we have
\begin{align}
\MO{n}&=\Szi{n}n\frac{J^n}{{\caln}^2(J)\rho_n}=
{\caln}^{-2}\Syx{n}{1}{\infty}\frac{(2J)^nn\skob{\half}_n
\skob{\frac32}_n}{n!(1)_n}= \nonumber\\
&=\frac{2J}{N^2}\Szi{n}\frac{(2J)^n\frac34\skob{\frac32}_n\skob{\frac52}_n
}{n!(2)_n}=
\frac{3J}{2}\,\,\frac{\Fto{\frac32}{\frac52}{2}{2J}}{\aFto}.
\label{Le.58}
\end{align}
Using the formulas (7.3.2 (50)) and (7.3.2 (217))from \cite{Prud},
we obtain the expression for $\MO{n}$ in terms of an elliptic
integrals
\begin{equation}\label{Le.59}
\MO{n}=\half\,\frac{J}{1-2J}\,
\frac{2\kbf(\sqrt{2J})-(1+2J)\dbf(\sqrt{2J})}{\ebf(\sqrt{2J})},
\end{equation}
where
\begin{align}
\ebf(k)&=\Ixy{0}{\frac{\pi}{2}}{\sqrt{1-k^2\sin^2t}}{t}
\qquad\qquad - \text{a full elliptic integral of a 2-nd kind},
\label{Le.60} \\
\dbf(k)&=\Ixy{0}{\frac{\pi}{2}}{\frac{\sin^2t}{\sqrt{1-k^2\sin^2t}}}{t}
=D(\frac{\pi}{2},k)\qquad - \text{a full elliptic integral},
\label{Le.61} \\
\kbf(k)&\!=\!\Ixy{0}\!{\frac{\pi}{2}}{\frac{1}{\sqrt{1-k^2\sin^2t}}}{t}
\!=\!F(\frac{\pi}{2},k) - \text{a full elliptic integral of a 1-st
kind}. \label{Le.62}
\end{align}

The similar evaluations give us that
\begin{align}
\MO{n^2}&=\Szi{n}n^2\frac{J^n}{N(J)^2\rho_n}=
N^{-2}\Syx{n}{1}{\infty}\frac{(2J)^n\half\skob{\frac32}^{n-1}\frac32
\skob{\frac52}^{n-1}}{(n-1)!(1)_{n-1}}=\nonumber\\
&=\frac32\frac{J}{N^2}\Szi{n}\frac{(2J)^n}{n!}
\frac{\skob{\frac32}_n\,\skob{\frac52}_n}{\skob{1}_n}= \frac32\,
J\,\frac{\Fto{\frac32}{\frac52}{1}{2J}}{\aFto}. \label{Le.63}
\end{align}
In terms of  elliptic integrals we can rewrite \reff{Le.63}
\begin{equation}\label{Le.64}
\MO{n^2}= \half\frac{J}{(1-2J)^2}
\frac{(3+10J)\kbf(\sqrt{2J})-2J(7+2J)\dbf(\sqrt{2J})}{\ebf(\sqrt{2J})}.
\end{equation}

The computation of the variance gives the following relation
\begin{equation}\label{Le.65}
\Delta n=\frac{\sqrt{\frac32J}}{\aFto}
\bskob{\bFto\aFto-\frac32J\skob{\dFto}^2}^{\half}.
\end{equation}
Then we obtain for a Mandel parameter $Q=\frac{(\Delta
n)^2}{\MO{n}}-1 $  a relation
\begin{equation}
Q=\frac{\bFto}{\dFto}-\frac32J\, \frac{\dFto }{\aFto}-1.
\end{equation}

\section{Conclusion.}
In the present note we have defined a new type of an oscillator
for which the Legendre polynomials play the same role as the
Hermite polynomials play for standard boson oscillator.Solving the
appropriate classical moment problem,we defined two sets of
coherent states - as eigenvectors of an annihilation operator and
temporary stable coherent state of the Klauder - Gazeau type. In
our next work which is in the closing stage , we shall define the
Perelomov type coherent states for the Legendre oscillator, as
once more concrete example of general study of connections between
orthogonal polynomials and coherent states. Similar
oscillator-like systems can be defined for others orthogonal
polynomials\cite{borz} (including $q$-deformed ones). For these
systems one can also define the corresponding systems of coherent
states (the rather general construction we will to describe in our
following work).

\end {document}